\documentclass[a4paper]{article}
\usepackage{latexsym}
\usepackage{amsmath}
\usepackage{amssymb}

\newtheorem{Theorem} {Theorem} [section]

\newtheorem{Conjecture} [Theorem] {Conjecture}
\newcommand{\Proof}{ \noindent{\bf Proof.}\quad }
\newcommand{\qed}{\hfill$\Box$\medskip}
\newcommand{\Ff}{{\mathbb F}}
\newcommand{\Qq}{{\mathbb Q}}

\author{Aart Blokhuis, Andries Brouwer, Benne de Weger \\
        Eindhoven University of Technology \\
        \texttt{a.blokhuis@tue.nl, aeb@cwi.nl, b.m.m.d.weger@tue.nl}}
\title{Binomial collisions and near collisions}
\date{\today}

\begin{document}
\maketitle

\begin{abstract}
We describe efficient algorithms to search for cases in which binomial coefficients
%BdW vervangen:
%are equal or almost equal, and give a conjecturally complete list of all cases
%where two binomial coefficients differ by 1.
are equal or almost equal, give a conjecturally complete list of all cases
where two binomial coefficients differ by 1, and give some identities
for binomial coefficients that seem to be new.
%BdW einde
\end{abstract}

\section{Introduction}
Let us call a quadruple $(n,k,m,l)$ with $2 \le k \le n/2$
and $2 \le l \le m/2$ a (binomial) {\em collision} when
$k < l$ and $\binom{n}{k} = \binom{m}{l}$,
and a {\em near collision} when $\binom{m}{l} - \binom{n}{k} = d > 0$
with $\binom{m}{l} \ge d^3$.
%BdW toegevoegd:
The exponent $ 3 $ is somewhat arbitrary. Maybe $ 5 $ is a more natural exponent, 
see the end of this paper.
%BdW einde

Collisions have been studied by many authors. Some references will
be given below. In this note we report on computer searches for
collisions and near collisions, and give seven infinite families
of near collisions.

\section{Collisions}
We list the known collisions. There are the double collision
\[
\binom{78}{2} = \binom{15}{5} = \binom{14}{6} = 3003,
\]
%BdW vervangen:
%six further sporadic examples given in the table below,
%and a miraculous infinite family.
six further sporadic examples given in the table below:
%BdW einde

% \medskip
\begin{center}
\begin{tabular}{ccccc}
$n$ & $k$ & $m$ & $l$ & $\binom{m}{l} = \binom{n}{k}$\\[0.5mm]
\hline
16 & 2 & 10 & 3 & 120 \\
21 & 2 & 10 & 4 & 210 \\
56 & 2 & 22 & 3 & 1540 \\
120 & 2 & 36 & 3 & 7140 \\
153 & 2 & 19 & 5 & 11628 \\
221 & 2 & 17 & 8 & 24310
\end{tabular}
\end{center}
%BdW vervangen
%
%The infinite family is given by
and a miraculous infinite family given by
%BdW einde
\[
\binom{F_{2i+2} F_{2i+3}}{F_{2i} F_{2i+3}} =
\binom{F_{2i+2} F_{2i+3}-1}{F_{2i} F_{2i+3}+1}
\quad \textrm{ for } \quad i = 1, 2, \ldots ,
\]
where $F_i$ is the $i$th Fibonacci number (defined by $F_0 = 0$, $F_1 = 1$,
$F_{i+1} = F_i + F_{i-1}$ for $i \ge 1$).
The infinite family is due to Lind \cite{L}, and was rediscovered by
several others such as Singmaster \cite{S2} and Tovey \cite{T}.
% Note that one may write
% \[
% F_{2i+2} F_{2i+3} = \dfrac{1}{5} \left( L_{4i+5} - 1 \right), \quad
% F_{2i} F_{2i+3} = \dfrac{1}{5} \left( L_{4i+3} - 4 \right),
% \]
% where $L_i$ is the $i$th Lucas number (defined by $L_0 = 2$, $L_1 = 1$,
% $L_{i+1} = L_i + L_{i-1}$ for $i \ge 1$).
%
Examples are
\[
\binom{15}{5} = \binom{14}{6},~~
\binom{104}{39} = \binom{103}{40},~~
\binom{714}{272} = \binom{713}{273},~~
\binom{4895}{1869} = \binom{4894}{1870}.
\]

\medskip
Twenty years ago one of us conjectured

\begin{Conjecture} {\rm (\cite{dW2})}
There are no other collisions than those given above.
\end{Conjecture}

The current status is as follows.

\begin{Theorem}
There are no unknown collisions in the following cases:
\begin{itemize}
\item $ (k,l) = (2,3), (2,4), (2,5), (2,6), (2,8), (3,4), (3,6), (4,6), (4,8) $,
\item $ (m,l) = (n-1,k+1), (n-1,k+2), (n-2,k+1) $,
\item $ n \leq 10^6 $,
\item $ \binom{n}{k} \leq 10^{60} $.
\end{itemize}
\end{Theorem}
\Proof
The first two parts can be found in the literature.

The case $ (k,l) = (2,3) $ was settled in \cite{A}.
The case $ (k,l) = (2,4) $ was settled in \cite{P}, and also in \cite{dW1}.
The case $ (k,l) = (2,5) $ was settled in \cite{BMSST}.
The cases $ (k,l) = (2,6), (2,8), (3,6), (4,6), (4,8) $ were settled in \cite{SdW1}.
The case $ (k,l) = (3,4) $ was settled in Mordell \cite{M} (actually, he solved
an equivalent equation and seems not to have noted the relation to binomial coefficients).

The case $ (m,l) = (n-1,k+1) $ was settled in \cite{T}
(and yields the infinite family).
The cases $ (m,l) = (n-1,k+2), (n-2,k+1) $ were settled in \cite{SdW2}.

The last two parts are the results of computer searches we report on in this paper.
Some details are given below.
\qed

Earlier computer searches handled $n \le 10^3$ and $\binom{n}{k} \le 10^{30}$ (\cite{dW2}).
In the literature one also finds finiteness results (\cite{F}, \cite{B}),
and results on the number of times an integer may occur as binomial coefficient
(\cite{S1}, \cite{AEH}, \cite{K}).

\subsection{Settling $n \le 10^6$}
In order to find all collisions with $n \le N$ for some fixed $N$,
generate a list of all values $ \binom{n}{k} $
with $ 2 \le k \le n/2 $ and $ n \le N $.
Sort it, and compare successive elements to find duplicates.

Now the list has length about $\frac14 N^2$, and probably does not fit into memory.
One approach is to split the list into parts, e.g.\ into the sublists
consisting of all binomial coefficients between $ 10^{e-1} $ and $ 10^e $,
for all relevant $ e $. We tried this in Mathematica and did $N = 34000$
in 23h30m on a 2.6 GHz Intel i7.

A different approach is to have a table and a priority queue, both of size $N$.
Both contain the same elements. Initially both contain the numbers
$\binom{n+2}{2}$ for $n < N$. The priority queue is kept sorted.
At each step the top two elements are compared for equality.
Afterwards the top element is discarded.
When $\binom{n+k}{k}$ is discarded, the new value $\binom{n+k+1}{k+1}$
is added, unless $k \ge n$. The new value needed is computed
from the old one via $\binom{n+k+1}{k+1} = \binom{n+k}{k} + \binom{n+k}{k+1}$.
Note that the value $\binom{n+k}{k+1}$ is present in the table at index $n-1$
at the moment it is needed.

Computation time for the algorithms as described is cubic in $N$
if the precise value of the binomial coefficients is computed,
since not only the length of the list grows, but also the size of the numbers.
% (For example, $\binom{1000000}{500000}$ has 301027 digits.)
Bounded precision suffices to ensure that (almost) collisions
are unlikely, and reduces the time needed to $O(N^2 \log N)$.
Almost collisions still occur (for example,
% $\binom{74372}{7254} = 5.588462942048745 \cdot 10^{10321}$,
% $\binom{88630}{6712} = 5.588462942048773 \cdot 10^{10321}$).
$\binom{102091}{12877} = 1.256839391954534 \cdot 10^{16800}$,
$\binom{200954}{9642} = 1.256839391954529 \cdot 10^{16800}$).
We used interval arithmetic to distinguish almost equal numbers,
and full exact multiple length arithmetic in the few cases where
the interval arithmetic did not suffice.
We tried this in C, with a custom data type (since the usual data types
do not handle large exponents, or are too slow), and did $N = 10^6$
in 56h14m on an old 2 GHz PC.

\subsection{Settling $\binom{n}{k} \leq 10^{60}$}
In order to find all collisions with $\binom{n}{k} \le M$ we handle
each relevant pair $(k,l)$ separately. Let $l_{\rm max}$ be the largest
$l$ such that $\binom{2l}{l} \le M$.
As we saw, the pairs $(k,l)$ with $k < l \le 4$ have been done already,
so it suffices to handle $5 \le l \le l_{\rm max}$, and for each
$l$ the values of $k$ with $2 \le k \le l-1$, with $k \ge 3$ if $l = 5$.

Given a pair $(k,l)$, let $m_{\rm max}$ be the largest $m$ with $\binom{m}{l} \le M$.
Make a list of all $m$ with $2l \le m \le m_{\rm max}$,
and discard the $m$ for which $\binom{m}{l}$ cannot be of the form $\binom{n}{k}$.
What is left are possible collisions, and in practice only actual collisions are left.

The discarding is done via a sieving process.
The function $f(n) = \binom{n}{k}$ is a polynomial of degree $k$
in the variable $n$, with rational coefficients.
The denominators of the coefficients have only prime factors $\leq k$.
For any prime $p > k$ the function $f$ induces a polynomial map from $\Ff_p$ to itself.
Let $A(k,p)$ be the size of the image. Experience shows that
$A(k,p) \approx (1-e^{-1})p$ when $k$ is odd, and $A(k,p) \approx (1-e^{-1/2})p$
when $k$ is even. See below for more remarks on this function $A(k,p)$.
Since $1-e^{-1} = 0.63...$ and $1-e^{-1/2} = 0.39...$, a significant fraction
of all residues mod $p$ cannot be of the form $\binom{n}{k}$.
Now pick $p > l > k$, and look at $\binom{a}{l}$ for $0 \le a \le p-1$.
Whenever $\binom{a}{l}$ is not in the mod $p$ image of $f$, discard
all $\binom{m}{l}$ with $m \equiv a$ (mod $p$) from the list.

Repeating this sieve action for all primes less than 500 (stopping earlier
when the list has become empty) we found all collisions up to $M = 10^{60}$.
The largest prime needed was $p = 401$. This took about 375 CPU hours total
on a few old 2 GHz machines. For large $l$ the upper bound $m_{\rm max}$
is small, and sieving is very quick. (In fact for $l \ge 10$ we sieved up to $10^{100}$.)
The main part of the work are the pairs $(k,l) = (3,5),(4,5)$, where the list
has length roughly $M^{1/5}$.

\subsubsection*{On $A(k,p)$}
There is a lot of literature on the size of the image of a polynomial on $\Ff_p$.
For $k = 3$ and $k = 4$ the value of $A(k,p)$ was found by
Daublebsky~v.~Sterneck~\cite{DvS}. One has $A(3,p) = (2p\pm1)/3$
when $p \equiv \pm1$ (mod 6), and
$A(4,p) = (3p+4+\chi(-1)+2\chi(5)-2\chi(10))/8$ for $p > 5$,
where $\chi$ is the quadratic character.
Birch \& Swinnerton-Dyer \cite{BSD} showed that `general' polynomials
of degree $k$ on $\Ff_p$ have an image of size $a_k p + O(\sqrt{p})$
where $a_k = \sum_{i=1}^k (-1)^{i-1} \frac{1}{i!}$.
We conjecture in our situation that the value
$A_k = \lim_{p\rightarrow\infty} \frac{A(k,p)}{p}$ exists,
and equals $A_k = a_k$ for odd $k$,
and $A_k = \sum_{i=1}^{k/2} (-1)^{i-1} \frac{1}{2^ii!}$ for even $k$.
This is true for $k \le 5$.
Note that for even $k$ there is the symmetry $f(x) = f(k+1-x)$ explaining
the smaller image size.

\section{Near collisions}

\subsection{Difference 1}
We know about the following examples with $d=1$:

\begin{center}
\begin{tabular}{ccccc}
$n$ & $k$ & $m$ & $l$ & $\binom{m}{l} = \binom{n}{k}+1$\\
\hline
    6 &  3 &       7 & 2 & 21 \\
    7 &  3 &       9 & 2 & 36 \\
   11 &  2 &       8 & 3 & 56 \\
   10 &  5 &      23 & 2 & 253 \\
   12 &  4 &      32 & 2 & 496 \\
   16 &  3 &      34 & 2 & 561 \\
   60 &  2 &      23 & 3 & 1771 \\
   27 &  3 &      77 & 2 & 2926 \\
   29 &  3 &      86 & 2 & 3655 \\
   34 &  3 &      21 & 4 & 5985 \\
   22 &  5 &     230 & 2 & 26335 \\
  260 &  3 &    2407 & 2 & 2895621 \\
   93 &  4 &    2417 & 2 & 2919736 \\
   62 &  5 &    3598 & 2 & 6471003 \\
   28 & 11 &    6554 & 2 & 21474181 \\
  665 &  3 &    9879 & 2 & 48792381 \\
  135 &  5 &   26333 & 2 & 346700278 \\
  139 &  5 &   28358 & 2 & 402073903 \\
19630 &  3 & 1587767 & 2 & 1260501229261 \\
160403633 & 2 & 425779 & 3 & 12864662659597529
\end{tabular}
\end{center}

The above table is complete for the cases
$(k,l),\,(l,k) = (2,3)$, $(2,4)$, $(2,6)$, $(3,4)$, $(4,6)$, $(4,8)$
and $ (k,l) = (2,8) $
(as one sees by finding all integral points on the corresponding elliptic
curves), and for $\binom{n}{k} \le 10^{30}$.
% All entries except the last were already given in \cite{mjqxxxx}.
% BdW vervangen
%We conjecture that this is the full list.
%More generally we conjecture that, given a fixed difference $d$,
%the number of near-collisions with difference $d$ is finite.
We conjecture the following
\begin{Conjecture} 
There are no other near collisions with difference $ 1 $ than those given above.
\end{Conjecture}
and, more generally,
\begin{Conjecture}
Given a fixed difference $d$, the number of near collisions with difference $d$ is finite.\end{Conjecture}
The latter conjecture can be backed by standard heuristic arguments.
%BdW einde
The infinite family of collisions seems like a miracle.
% For $d=2$ only $k=2 l=3 i=104, mm=182104$

%The cases mentioned above
% and (2,8), (4,8)
%correspond to (double covers of) elliptic curves in cubic Weierstrass form
%or in quartic form, and can be solved completely using
%e.g.\ Sage \cite{Sage} or Magma \cite{Magma}.
%See \cite{SdW1}, Table 1, for the transformations from the
%binomial equations to the elliptic equations.

The cases $(k,l),\,(l,k) = (2,3)$, $(2,4)$, $(2,6)$, $(2,8)$, $(3,4)$, 
$(3,6)$, $(4,6)$, $(4,8)$ correspond to integral points on (double 
covers of) elliptic curves, that can in principle be solved by the 
methods of \cite{SdW1}, \cite{SdW2}. All except $(3,6)$ are curves in 
Weierstrass or quartic form, and can in principle be solved completely 
using e.g.\ Sage \cite{Sage} or Magma \cite{Magma}. We succeeded in
doing so using Magma for all except $(k,l) = (8,2)$.
See \cite{SdW1}, Table 1, for the transformations from the
binomial equations to the elliptic equations.

\subsection{Infinite families}
When $d$ is not fixed, there are a few infinite families:
\begin{align}
& \binom{12x^2-12x+3}{3} + \binom{x}{2} = \binom{24x^3-36x^2+15x-1}{2} \\
& \binom{12x^2-12x+5}{3} + \binom{x}{2} = \binom{24x^3-36x^2+21x-4}{2} \\
& \binom{60x^2-60x+15}{5} + \binom{x}{2} = \binom{a}{2}
\intertext{\rm where $a = 3600x^5-9000x^4+8700x^3-4050x^2+905x-77$,}
& \binom{60x^2-60x+19}{5} + \binom{x}{2} = \binom{a}{2}
\intertext{\rm where $a = 3600x^5-9000x^4+9300x^3-4950x^2+1355x-152$,}
& \binom{240x^2-240x+62}{5} + \binom{3x-1}{2} = \binom{a}{2}
%
% \intertext does not handle \\, but \newline works
%
\intertext{\rm where
$a = 115200x^5-288000x^4+288000x^3-144000x^2+35995x-3597$,}
& \binom{11340x^2+11340x+2835}{9} + \binom{y}{2} = \binom{a}{2}
\intertext{where $y = 22680x^3+34020x^2+17001x+2831$ and
$a = 4134207084840000x^9+18603931881780000x^8+37201301530092000x^7+
43386206573682000x^6+$\newline
$32522432635935900x^5+16249739546454750x^4+5411800833695550x^3+$\newline
$1158443736409575x^2+144626588131776x+8023467184451$,}
& \binom{11340x^2+11340x+2843}{9} + \binom{y}{2} = \binom{a}{2}
\end{align}
where $y = 22680x^3+34020x^2+17019x+2840$ and
$a = 4134207084840000x^9+18603931881780000x^8+37214425997028000x^7+
43432142207958000x^6+$\\
$32591336087349900x^5+16307159089299750x^4+5440510606648950x^3+$\\
$1167056670132675x^2+146062077851076x+8126002273751$.

\medskip
How does one find such identities?
In order to get $\binom{b}{3} + \binom{x}{2} = \binom{a}{2}$,
where $x$ is small, one needs $\frac{1}{3}b(b-1)(b-2) = (a-x)(a+x-1)$,
a product of two nearly equal numbers. If $b = 3e^2$, then
$\frac{1}{3}b(b-1)(b-2) = (e(b-2))(e(b-1))$ and we can take
$a-x = e(b-2)$, $a+x-1 = e(b-1)$ and find $e = 2x-1$, $b = 3(2x-1)^2$,
the first family. The other families arise in a similar way.

% Similarly, consider $\binom{b+4}{9}=\binom{a}{2}-\binom{x}{2}$ and write
% the left hand side as as $c_1c_2c_3$, where $c_1=64(b-4)/9!=(2e+1)^2/2$,
% $8c_2=(b-2)(b-1)(b+1)(b+4)=b^4+2b^3-9b^2-2b+8$ and
% $8c_3=(b-3)b(b+2)(b+3)=b^4+2b^3-9b^2-18b$.
% We now take $a+x-1=(2e+1)c_2$ and $a-x=(2e+1)c_3$, and find
% $x=(2e+1)b+e+1$ and $a=(2e+1)(b^4+2b^3-9b^2-10b)/8+e+1$,
% where $b-4 = 3^4.5.7(2e+1)^2$.
% % e=0: b=2839, x=2840, a=8126002273751
% % since b is odd, (b^4+2b^3-9b^2-10b) is divisible by 8 - factor b^2-1 mod 8.
% We found
% \begin{align}
% & \binom{2835(2e+1)^2+8}{9} + \binom{x}{2} = \binom{a}{2}
% \end{align}
% where the difference $\binom{x}{2}$ is less than the cube root
% of $\binom{a}{2}$.
%
% Now the mirror image:
% $c_1=64(b+4)/9!=(2e+1)^2/2$
% $8c_2=(b+2)(b+1)(b-1)(b-4)=b^4-2b^3-9b^2+2b+8$
% $8c_3=(b+3)b(b-2)(b-3)=b^4-2b^3-9b^2+18b$
% $a+x-1=(2e+1)c_3$ and $a-x=(2e+1)c_2$ gives
% $x=(2e+1)b-e$ and $a=(2e+1)(b^4-2b^3-9b^2+10b)/8+e+1$
% where $b+4 = 3^4.5.7(2e+1)^2$.
% % e=0: b=2831, x=2831, a=8023467184451

Are there many such identities?
Let us say that the quality of an identity
$\binom{n(x)}{k} + d(x) = \binom{m(x)}{l}$
is the degree of $x$ in $\binom{n(x)}{k}$ and $\binom{m(x)}{l}$
divided by that in $d(x)$. Then our identities (1)-(7)
have qualities 3, 3, 5, 5, 5, 3, 3. These are the only identities
of quality at least 3 that we know of. Maybe there are no others.
%BdW toegevoegd:
Maybe there is a number $ \alpha $, supposedly $ \leq 3 $, such that there are
only finitely many identities of quality at least $ \alpha $. 

It follows from the existence of the identities (1)-(7) that there are 
infinitely many near collisions, even with $ \binom{m}{l} \geq d^5 $.
Maybe there is a number $ \beta $, certainly $ \beta > 5 $, such that there are
only finitely many near collisions with $ \binom{m}{l} \geq d^{\beta} $.
%BdW einde

\end{document}